\documentclass[twocolumn]{autart}    

\usepackage{graphicx,amsmath}          

\usepackage[dvipsnames]{xcolor}

\newtheorem{theorem}{Theorem}[section]

\newtheorem{e-proposition}[theorem]{Proposition}
\newtheorem{corollary}[theorem]{Corollary}
\newtheorem{e-definition}[theorem]{Definition\rm}
\newtheorem{remark}{\it Remark\/}

\newtheorem{theoreme}{Th\'eor\`eme}[section]

\newtheorem{proposition}[theoreme]{Proposition}

\renewcommand{\theequation}{\arabic{equation}}
\begin{document}

\def\blue#1{\textcolor{blue}{#1}}
\def\red#1{\textcolor{red!70!black}{#1}}

\let\mib=\boldsymbol
\let\cal=\mathcal
\newcommand{\argmax}{\arg\!\max}

\def\mA{{\bf A}} 
\def\mB{{\bf B}} 
\def\mG{{\bf G}} 
\def\mI{{\bf I}}
\def\mmu{{\mib \mu}}
\def\mnu{{\mib \nu}}
\def\mnul{{\mib 0}}
\def\mx{{\bf x}}

\def\vc{{\sf c}}
\def\vnul{{\sf 0}}
\def\vu{{\sf u}} 
\def\vv{{\sf v}} 
\def\vx{{\sf x}} 
\def\vy{{\sf y}} 
\def\vz{{\sf z}} 

\def\card{{\rm card}}
\def\dist{{\rm dist}}
\def\Dup#1#2{\langle#1,#2\rangle}
\def\dv{{\sf div\thinspace}}
\def\eps{\varepsilon} 
\def\HH#1{{{\rm H}^{#1}}}
\def\HHj{{{\rm H}^1}}
\def\lA{{\cal A}}
\def\lF{{\cal F}}
\def\lK{{\cal K}}
\def\lN{{\cal N}}
\def\lO{{\cal O}}
\def\LdTT{{{\rm L}^2(0,T)}} 
\def\LLd{{{\rm L}^2}}
\def\LLb{{{\rm L}^{\infty}}}
\def\nor#1#2{{\| #1 \|}_{#2}}
\def\Nor#1{\| #1 \|}
\def\oi#1#2{\langle#1,#2\rangle}
\def\ozi#1#2{\langle#1,#2]}
\def\ph{\varphi} 
\def\Rpl{{{\bf R}^{+}}}
\def\span{{\rm span}}
\def\svaki#1{(\forall\,#1)}
\def\Svaki#1{\left(\forall\,#1\right)}
\def\zi#1#2{[#1,#2]}

\def\zz{\mathbb{Z}} 
\def\rr{{\bf R}}
\def\nn{\mathbb{N}}
\def\cc{\mathbb{C}}
\def\C{{\bf C}}
\def\N{{\bf N}}
\def\Z{\bf Z}
\def\R{{\bf R}}
\def\G{\Gamma}
\def\Th{\Theta}
\def\L{\Lambda}
\def\Si{\Sigma}
\def\F{\Phi}
\def\O{\Omega}
\begin{frontmatter}


  \title{Greedy  controllability of finite dimensional linear systems
    \thanksref{footnoteinfo}} 

\thanks[footnoteinfo]{
This paper has been partially developed while the first author was visiting the Basque Center for Applied Mathematics (Bilbao, Spain) within the MTM2011-29306 Grant of the MINECO. 
  The first author was also supported  by Croatian Science
Foundation under the project  WeConMApp/HRZZ-9780, as well as  by University of Dubrovnik trough the grant 11-18-10-2 and through the Erasmus+ programme.
\\
The second author was partially supported
by the Advanced Grant NUMERIWAVES/FP7-246775 of
the European Research Council Executive Agency,
 the FA9550-15-1-0027 of AFOSR, the
MTM2011-29306 and MTM2014-52347 
Grants of the MINECO, and a Humboldt Award at
the University of Erlangen-N\"urnberg. 
 }
\author[UNIDU]{Martin Lazar}
\author[UAM]{Enrique Zuazua}
\ead{martin.lazar@unidu.hr}
\address[UNIDU]{University of Dubrovnik, 
\'Cira Cari\' ca 4, 20 000 Dub\-rov\-nik, Croatia }

\ead{enrique.zuazua@uam.es}
\address[UAM]{Departamento de Matem\'aticas,
Universidad Aut\' onoma de Madrid,
Cantoblanco,
28049 Madrid, Spain}


\begin{abstract}                          
We analyse the problem of controllability for  parameter dependent linear finite-dimensional systems. The goal is  to identify the most distinguished realisations of those parameters so to better describe or approximate the whole range of controls. We adapt  recent results on greedy and weak greedy algorithms for parameter dependent PDEs or, more generally, abstract equations in Banach spaces.  Our results lead to optimal approximation procedures that, in particular, perform better than simply sampling the parameter-space to compute the controls for each of the parameter values. We apply these results for the approximate control of finite-difference approximations of the heat and the wave equation. The numerical experiments  confirm the efficiency of the methods and show that the number of weak-greedy samplings that are required is particularly low when dealing with  heat-like equations,  because of the intrinsic dissipativity that the model introduces for high frequencies.
\end{abstract}

\end{frontmatter}

\abovedisplayskip=4pt
\belowdisplayskip=4pt


\section{Introduction and problem formulation}
We analyse the problem of controllability for linear  finite-dimensional systems submitted to parametrised perturbations,  depending on unknown parameters in a deterministic manner. 

In previous works we have analysed the property of averaged control looking for a control, independent of the values of these parameters,  designed to perform well, in an averaged  sense (\cite{Z13},  \cite{LZ}).

Here we analyse the complementary issue of determining the most relevant values of the unknown parameters so to provide the best possible approximation of the set of parameter dependent controls. Our analysis is based on previous work on reduced modelling and (weak) greedy algorithms for parameter dependent PDEs and abstract equations in Banach spaces (\cite{Buffa12}, \cite{CDV}), which we adapt to the present context.

  The problem is relevant in applications, as in practice  the models under consideration are often not
completely determined, submitted to unknown or uncertain parameters, either of
deterministic or of stochastic nature. It is therefore essential to develop robust analytical and computational methods, not only allowing to control a given
model, but also to deal with parameter-dependent families of systems in a
stable and computationally efficient way.

Both reduced modelling and the control theory have
experienced successful real-life implementations (we refer to the book \cite{Le} for a series of such interactions with the European industry). The merge of this two
theories will allow of variety of applications in all fields involving problems modelled by
parameter dependent systems (fluid dynamics, aeronautics,
meteorology, economics, ...).

Although the  greedy control  is applicable to more general control problems and systems, here we  concentrate on controllability issues and, to better illustrate the main ideas of the new approach, we focus on linear finite-dimensional systems of parameter dependent ODEs.  Infinite-dimensional systems, as a first attempt to later consider PDE models,  are discussed separately in Section \ref{infinite}, as well as in the Conclusion section.

Consider the finite dimensional linear control system
\begin{equation}\label{eq1F-d}
\left\{\begin{array}{ll} \frac{d}{dt}x(t, \nu)=\mA(\nu)x(t, \nu)+\mB(\nu)u(t, \nu), \, 0<t<T,\\
x(0)=x^0.
\end{array}\right.
\end{equation}

In (\ref{eq1F-d}) the (column) vector valued function
$x(t, \nu)=\big(x_1(t, \nu),\ldots,x_N(t, \nu)\big)\in {\bf R}^N$ is the state
of the system
at time $t$ governed by dynamics determined by the parameter $\nu \in \lN \subset {\bf R}^d$, $d\ge 1$, $\lN$ being a compact set,
$\mA(\nu)$ is a $N \times N-$matrix governing its free dynamics and $u=u(t, \nu)$ is a $M$-component control vector in ${\bf R}^M$, $M \le N$, entering and acting on the system through the control operator $\mB(\nu)$, a $N \times M$ parameter dependent matrix. In the sequel, to simplify the notation, $d/dt$ will be simply denoted by $'$.

The matrices $\mA$ and $\mB$ are  assumed
to be Lipschitz continuous with respect to the parameter $\nu$. However, some of our analytical results (Section \ref{infinite}) will additionally require analytic dependence conditions on $\nu$.

Here, to simplify the presentation, we have assumed   the initial datum $x^0 \in  {\bf R}^N$ to be controlled, to be independent of the parameter $\nu$. Despite of this, the matrices $\mA$ and $\mB$ being $\nu$-dependent, both the control and the solution will depend on $\nu$. 
Similar arguments allow to handle the case when $x^0$ also depends on the parameter $\nu$, which will be discussed separately. 

We address the controllability of this system whose initial datum  $x^0$ is given, known and fully determined. 
We assume that the system under consideration is controllable for all values of $\nu$. This can be ensured to hold, for instance, assuming that the controllability condition is satisfied for some specific realisation $\nu_0$ and that the variations of $\mA(\nu)$ and $\mB(\nu)$ with respect to $\nu$ are small enough.

In these circumstances, for each value of $\nu$ there is a control of minimal $[L^2(0, T)]^M$-norm, $u(t, \nu)$. This defines a map, $\nu \in \lN \to [L^2(0, T)]^M$, whose regularity is determined by that of the matrices entering in the system, $\mA(\nu)$ and $\mB(\nu)$.

Here we are interested on the problem of determining the optimal selection of a finite number of realisations  of the parameter $\nu$ so that all controls, for all possible values of $\nu$, are optimally approximated.

More precisely, the problem can be formulated as follows.

{\bf Problem 1}
{\it Given a control time
$T>0$ and arbitrary initial data $x^0$ and final target $x^1\in{\bf R}^N$,  we consider the set of controls of minimal $[L^2(0, T)]^M$-norm, $u(t, \lN)$, corresponding to all possible values $\nu \in \lN$ of the parameter satisfying the controllability condition:
\begin{equation}\label{eq2F-d}
 x(T, \nu)  =x^1.
\end{equation}
This set of controls is compact   in $[L^2(0, T)]^M$. 

Given $\varepsilon > 0$ we aim at determining   a family of parameters $\nu_1,...,  \nu_{n}$ in $\lN$,  whose cardinal $n$ depends on $\eps$,  so that 
the corresponding controls,  denoted by $u_1, ..., u_{n}$,  are such that for every $\nu \in \lN$ there exists $u^\star_\nu \in {\rm span}\{u_1,..., u_n\}$
steering the system \eqref{eq1F-d} in time $T$ within the $\eps$ distance from the target $x^1$, i. e. such that
\begin{equation}\label{eq2F-dap}
 ||x(T, \nu)  - x^1||\le \eps.
\end{equation}

}

Here and in the sequel, in order to simplify the notation, we  denote by $u_\nu$ the control $u(t, \nu)$, and similarly we use the simplified notation $\mA_\nu, \mB_\nu, x_\nu$.

Note that, in practice, the controllability condition (\ref{eq2F-d}) is relaxed to the approximate one (\ref{eq2F-dap}).
This is so since, in practical applications, when performing numerical approximations, one is interested in achieving the targets within a given error. This fact is also intrinsic to the methods we employ and develop in this paper, and that can only yield optimal procedures to compute approximations of the exact control, which turn out to be approximate controls in the sense of (\ref{eq2F-dap}).

  This problem is motivated by the practical issue of avoiding the
construction of a control function $u_\nu$ for each new parameter value $\nu$ which,  for large systems, although theoretically  feasible by the uniform controllability assumption,  would be computationally  expensive. By the contrary, the methods we develop try to exploit the advantages that a suitable choice of the most representative values of $\nu$ provides when computing rapidly the approximation of the control for any other value of $\nu$, ensuring that  the system is steered to the target   within the given error  (\ref{eq2F-dap}).
 
   Of course, the compactness of the parameter set $\lN$ and the  Lipschitz-dependence assumption with respect to $\nu$ make the goal to be feasible. It would suffice, for instance, to apply a {\it naive} approach, by taking a fine enough uniform  mesh on $\lN$ to achieve the goal.
However, our aim is to minimise the number of spanning controls $n$ and to derive the most efficient approximation. The {\it naive} approach is not suitable in this respect.

To achieve this goal we adapt to the present frame of finite-dimensional control, the theory developed in recent years based on greedy and weak-greedy algorithms for parameter dependent PDEs or abstract equations in Banach spaces,
which optimise the dimension of the approximating space, as well as the number of steps required for its construction.

The rest of this paper is organised as follows. In Section \ref{control-prel} we summarise the needed controllability results for  finite-dimensional systems and reformulate Problem 1 in terms  of the corresponding Gramian operator.
Section 3 is devoted to the review of  (weak) greedy algorithms, while their application to the control problem under consideration and its solution is provided in the subsequent section.

  The computational cost  of the greedy control approach is analysed in  Section 5. Section 6  contains a generalisation of the approach to infinite dimensional problems followed by  a convergence analysis of the greedy approximation errors with respect to the dimension of the approximating space.
Section 7 contains  numerical examples and experiments for finite-difference   discretisations of   1-D wave and heat problems. The paper is closed  pointing towards future development   lines of the greedy control approach.

\section{Preliminaries on finite dimensional control systems. Problem reformulation}
\label{control-prel}

In order to develop the analysis in this paper it is necessary to derive a convenient characterisation of the control of minimal norm $u_\nu$, as a function of the parameter $\nu$. This can be done in a straightforward manner in terms of the Gramian operator.
In this section we briefly summarise the most basic material on
finite-dimensional systems that will be used along this article
(we refer to \cite{MZ2,Z2} for  more details).

Consider the finite-dimensional system of dimension $N$:
\begin{equation}
\label{fd}
x'=\mA x+\mB u, \quad 0\le t \le T; \quad x(0)=x^0,
\end{equation}
where $x$ is the $N$-dimensional state and $u$ is the
$M$-dimensional control, with $M \le N$.

This corresponds to a specific realisation of the system above for a given choice of the parameter $\nu$. We omit however the presence of $\nu$ from the notation since we are now considering a generic linear finite-dimensional system.

Here $\mA$ is an $N\times N$ matrix with constant real coefficients
and $\mB$ is an $N\times M$ matrix. The matrix $\mA$ determines the
dynamics of the system and the matrix $\mB$ models the way $M$ controls
act on it.

In practice, it is desirable to control the $N$
components of the system with a low number of controls, the best possible case being the one of scalar controls: 
$M=1$.

Recall that system  (\ref{fd}) is said to be {\it controllable}  when
every initial datum $x^0\in \R^N$ can be driven to any final datum
$x^1$ in $\R^N$ in time $T$. This controllability property can be characterised 
by a
necessary and sufficient condition, which is of
purely algebraic nature, the so called {\it Kalman
condition}: {\rm System (\ref{fd})  is controllable if and only if
\begin{equation}
\label{rc}
{\rm{rank}}[\mB, \mA \mB, ..., \mA^{N-1}\mB]=N.
\end{equation} }

When this rank condition is fulfilled the system is controllable  for all $T>0$.

There is a direct proof of this result which uses the
representation of solutions of (\ref{fd}) by means of the
variations of constants formula. But for our purpose it is more convenient to use the point of view based on the dual problem of observability of the adjoint system that we discuss now.

Consider the {\it adjoint system}
\begin{equation}
  \label{afd}
  \left\{\begin{array}{ll}
  -\varphi'=\mA^*\varphi, \quad 0\le t \le T, \\
\varphi(T)=\varphi^0.
  \end{array}\right.
\end{equation}

\begin{proposition}
System (\ref{fd}) is controllable
in time $T$ if and only if the adjoint system (\ref{afd}) is {\it
observable}  in time $T$, i. e.  if  there exists a constant $C>0$
such that, for all solution $\varphi$ of (\ref{afd}),
\begin{equation}
\label{fdoi} 
|\varphi^0|^2 \le C \int_0^T |\mB^*\varphi|^2 dt.
\end{equation}
Both properties hold in all time $T$ if and only if the Kalman rank condition (\ref{rc}) is satisfied.

Furthermore, the control of minimal $[L^2(0, T)]^M$-norm can be built as a minimiser of a  quadratic functional $J: \R^N \to \R$:
\begin{equation}
  \label{funct_J}
J(\varphi^0)=\frac{1}{2}\int_0^T \!\!|\mB^*\varphi(t)|^2 dt - \Dup{x^1}{\varphi^0}+\Dup{x^0}{\varphi(0)}.
\end{equation}
More precisely, if\, $\tilde \varphi^0$ is a minimiser for
$J$,  then the control
\begin{equation}
 \label{control}
 u=\mB^*\tilde\varphi,
\end{equation}
where $\tilde\varphi$
is the corresponding solution of (\ref{afd}), is such that the state $x$ of (\ref{fd})
satisfies the control requirement $x(T)=x^1$.  

\end{proposition}

Note that the minimiser of $J: \R^N \to \R$ exists  since the functional $J$ is continuous, quadratic and coercive, in view of the observability inequality.

 This characterisation of controls ensuring controllability also yields explicit bounds on the controls. Indeed, since the functional $J \le 0$ at the minimiser, and in view of the observability inequality (\ref{fdoi}), it follows that
\begin{equation*}
||u|| \le 2\sqrt{C} [|x^0|^2 + |x^1|^2]^{1/2},
\end{equation*}
$C$ being the same constant as in (\ref{fdoi}). Therefore, we see that the square root of the observability constant is, up to a multiplicative factor, the norm of the control map associating to the data the control of minimal norm. 

Summarising, the control of minimal $[L^2(0, T)]^M$-norm is of the form $ u=\mB^*\tilde\varphi, $ where $\tilde\varphi= \exp(\mA^*(T-t)) \tilde \varphi^0$, $\tilde \varphi^0$ being the minimiser of $J$.

Let $\Lambda$ be the  quadratic form, known as  (controllability) Gramian,  associated to the pair $(\mA, \mB)$, i.e. 
$$
\Dup{\Lambda \varphi^0}{\psi^0}= \int_0^T \Dup{\mB^*\varphi }{ \mB^*\psi} dt,
$$
with $\varphi, \psi$ being solutions to \eqref{afd} with the data $\varphi^0$ and $\psi^0$, respectively. 
Then, under the rank condition, because of the observability inequality, this operator is coercive and symmetric and therefore invertible. Its corresponding matrix, which we denote the same, is given by the relation
\begin{equation}
\label{Gramian}
\L=\int_0^T e^{(T-t)\mA}\mB \mB^\ast e^{(T-t)\mA^\ast} dt\,.
\end{equation}

The minimiser $\tilde \varphi^0$ can be expressed as the solution to the linear system
\begin{equation*}
\label{phi^0}
\L   \tilde\ph^0 = x^1 -  e^{T\mA} x^0.
\end{equation*}
Hereby, the left hand side, up to the free dynamics component,  represents the solution of the control system (\ref{fd}) with the control given by \eqref{control}. As the solution is steered to the target $x^1$, the last relation follows.

In our context the adjoint system depends also on the parameter $\nu$:
\begin{equation}\label{eq-varp}
\left\{\begin{aligned}
-\varphi_\nu'(t)&=\mA_\nu^\ast (\nu) \varphi_\nu(t), \quad t\in(0,\,T),\\
\varphi_\nu(T)&=\varphi^0.
\end{aligned}\right.
\end{equation}
We assume that the system under consideration is controllable for all values of $\nu$. This can be ensured to hold assuming  the following uniform controllability condition
\begin{equation}
\label{G-bound}
\L_- \mI \leq  \L_\nu \leq \L_+ \mI,
\end{equation}
where $\L_{\pm}$ are positive constants, while $\L_\nu$ is the Gramian  of the system determined by $\mA_\nu, \mB_\nu$.

As we mentioned above, this assumption is fulfilled, in particular,  as soon as the system is controllable for some specific value of the parameter $\nu =\nu_0$, $\mA$ and $\mB$ depend continuously on $\nu$, and $\nu$ is close enough to $\nu_0$. But our discussion and presentation makes sense in the more general setting where (\ref{G-bound}) is fulfilled.

As we restrict the analysis to the set of minimal $[\LLd(0,T)]^M$-norm controls, each such control can be uniquely determined by the relation
\begin{equation}
\label{control_nu}
u_\nu = \mB_\nu^\ast e^{(T-t)\mA_\nu^\ast}  \ph^0_\nu,
\end{equation}
where $\ph^0_\nu$ is the unique  minimiser of a  quadratic functional $J_\nu: \R^N\to \R$ given by \eqref{funct_J}, with $\mB$ replaced by $\mB_\nu$, and $\ph$ by $\ph_\nu$, the solution to \eqref{eq-varp}.

As explained above, the minimiser $\ph^0_\nu$  can be equivalently  determined as the solution to the system 
\begin{equation}
\label{phi^0_par}
\L_\nu   \ph^0_\nu = x^1 -  e^{T\mA_\nu} x^0,
\end{equation}
where $\L
_\nu $  is the Gramian  associated to  $(\mA_\nu, \mB_\nu)$.

According to the uniform controllability condition \eqref{G-bound}, this defines a  mapping $\ph^0: \nu\in \lN \to \R^N $,
 whose smoothness is transferred from the mappings of $\mA_\nu$ and $\mB_\nu$ at all levels: Lipschitz, analytic etc.
Having assumed these maps are Lipschitz continuous and the parameter $\nu \in \lN$  varies on a compact set, the set of minimisers $\ph^0(\lN)$  constitutes also a compact  set in $\R^N$.

By using the 1-1 correspondence between the controls  $u_\nu$ and the associated minimisers $\ph^0_\nu$, the original problem formulation can be reduced to the following one.

{\bf Problem 2}
{\it Given a control time
  $T>0$, an arbitrary initial data $x^0$,  final target $x^1\in{\bf R}^N$,  and  $\varepsilon > 0$, and taking into account that the set of minimisers $\ph^0(\lN)$ corresponding to all possible values $\nu \in \lN$ of the parameter is compact in $\R^N$, 
 we aim at determining   a family of parameters $\nu_1,..., , \nu_n$ in $\lN$ so that 
the corresponding minimisers, that we denote by $\ph^0_1, ..., \ph^0_n$,  are such that for every $\nu \in \lN$ there exists $\ph^{0,\star}_\nu \in \span\{\ph^0_1,..., \ph^0_n\}$
such that the control $u_\nu^\star$ given by \eqref{control_nu}, with $\ph_\nu^0$ replaced by $\ph^{0,\star}_\nu$, steers the system \eqref{eq1F-d} to the state $x_\nu^\star(T)$ within the $\eps$ distance from the target $x^1$.  
}

In this formulation  the elements of the manifold we want to approximate are $N$-dimensional vectors (instead of $[\LLd(0, T)]^M$-functions in the first formulation), uniquely determined as solutions to linear systems \eqref{phi^0_par}. This enables to adapt the  (weak) greedy algorithms and reduced bases methods for parameter dependent problems, that we present in the following section.

 The approximate controls we obtain in this manner  do not really belong to the space $\span\{u_1, \ldots,  u_n\}$ spanned by controls associated to selected parameter values, since the selection is done at the level of $\ph^{0,\star}_\nu \in \span\{\ph^0_1,\ldots, \ph^0_n\}$, the control being simply the natural one corresponding to the choice of $\ph^{0,\star}_\nu $.

\section{Preliminaries on (weak) greedy algorithms}

In this section we present a brief introduction and main results of the linear approximation theory of parametric problems based on the  (weak) greedy algorithms, that we shall use in our application to controllability problems. A more exhaustive overview can be found in some  recent papers, e.g. \cite{CD15,DV15}.

The goal is to approximate a  compact set $\lK$ in a Banach space $X$ by a sequence of finite dimensional subspaces $V_n$ of dimension $n$. By increasing $n$ one improves the accuracy of the approximation.

Determining {\it offline} an approximation subspace within a given error normally implies a high computational effort. However, this calculation is performed only once, resulting in a good subspace from which one can easily and computationally cheaply construct {\it online}  approximations to every vector from  $\lK$. 

Vectors $x_i, i=1...n$ spanning the space $V_n$  are called {\it snapshots} of $\lK$. 

 The goal of  (weak) greedy algorithms is to construct a family of finite dimensional spaces $V_n\leq X$ that approximate the set $\lK$ in the best possible manner. 
The algorithm is structured  as follows. 

{\bf Weak greedy algorithm}

\vskip -3mm
\noindent\makebox[\linewidth]{\rule{0.5\textwidth}{0.4pt}}
Fix a constant $\gamma \in \ozi 01$. 

In the first step choose  $x_1 \in \lK$ such that 
\begin{equation}
\label{u_1}
\nor{x_1}X \geq \gamma  \max_{x \in \lK} \nor{x}X .
\end{equation}
At the general step, having found $x_1 ... x_n$, denote $V_n={\rm span} \{x_1, \ldots,x_n\} $ and
\begin{equation}
\label{greedy_error}
\sigma_n(\lK):= \max_{x \in \lK} {\rm dist}(x, V_n)\,.
\end{equation}
Choose the next element $x_{n+1}$ such that
\begin{equation}
\label{general-step}
{\rm dist}(x_{n+1}, V_n) \geq \gamma \sigma_n(\lK).
\end{equation}
The algorithm stops when $\sigma_n(\lK) $ becomes less than the given tolerance $\eps$.\\
\noindent\makebox[\linewidth]{\rule{0.5\textwidth}{0.4pt}}

The algorithm produces a finite dimensional space $V_n$ that approximates the set $\lK$ within the given tolerance $\eps$. The choice of a new element $x_n$ in each step is not unique, neither is the sequence of approximation rates $\sigma_n(\lK)$. But every such a chosen sequence decays at the same rate, which under certain assumptions given below, is close to the optimal one.
Thus the algorithm optimises  the number of steps required in order to satisfy the given tolerance, as well  as the dimension of the final space $V_n$.

The pure greedy algorithm corresponds to the case $\gamma=1$. As we shall see below, the relaxation of the pure greedy method ($\gamma=1$) to a weak greedy one ($\gamma \in \ozi 01$) will not significantly reduce the efficiency of the algorithm, making it, by the contrary, much easier for implementation.

When performing the (weak) greedy algorithm one has to chose the next element of the approximation space by exploring the distance \eqref{greedy_error} for all possible values $x\in \lK$. Such approach is faced with two crucial obstacles:
\begin{itemize}
\item[{\it i)}] the set $\lK$ in general consists of infinitely many vectors.
  \item[{\it ii)}] in practical implementations the set $\lK$ is often unknown (e.g. it represents the family of solutions to parameter dependent problems). 
\end{itemize}

The first problem  is bypassed by performing a search over some finite discrete subset of $\lK$. Here we use the fact that $\lK$, being a compact set, can be covered by a finite number of balls of an arbitrary small radius.

As to deal with the second one, instead of considering the exact distance appearing in \eqref{greedy_error}, one uses some {\it surrogate},  easier to compute.

In order to estimate  the efficiency of the weak greedy algorithm  we compare its approximation rates $\sigma_n(\lK)$ with the best possible ones.

The best choice of a approximating space $V_n$ is the one producing the smallest approximation error. This smallest error for a compact set $\lK$ is called the {\it Kolmogorov $n$-width of $\lK$}, and is defined as
$$
d_n(\lK):=\inf_{\dim Y=n} \sup_{x\in \lK} \, \inf_{y\in Y} \nor {x-y}X\,.
$$
It measures how well $\lK$ can 
be approximated by a subspace in $X$  of a fixed dimension $n$.

In the sequel we want to compare $\sigma_n(\lK)$ with the Kolmogorov width $d_n(\lK)$, which represents the best possible approximation of $\lK$ by a $n$ dimensional subspace of the referent Banach space $X$.

A precise estimate in that direction was provided by \cite{BCDDPW} in the Hilbert space setting,  and	subsequently improved and extended to the case of a general Banach space $X$  in \cite{DPW13}.

\begin{theorem}
  [\cite{DPW13}, Corollary 3.3]
  \label{greedy_rates}
  For the weak greedy algorithm with constant $\gamma$ in a Hilbert space $X$  we have the following:
If the compact set $\lK$ is such that, for some $\alpha>0$ and $C_0>0$
$$
    d_n(\lK)\leq C_0 n^{-\alpha}, \; n\in\N,
    $$
    then
    \begin{equation}
  \label{poly}
\sigma_n(\lK) \leq C_1 n^{-\alpha}, \; n\in\N,
    \end{equation}
    where $C_1:= \gamma^{-2} 2^{5\alpha+1} C_0$.
\end{theorem}

    This theorem implies that the weak greedy algorithms preserve  the polynomial  decay rates of the approximation errors, and the result will be   used  in the convergence analysis of our method in Section \ref{infinite} . A similar estimate also holds for  exponential  decays (cf. \cite{DPW13}).  In addition, it is also remarkable that
  the constant $\gamma$ effects the efficiency only up to a multiplicative constant, while leaving approximation rates unchanged.

\section{Greedy control and  the main result}

In this section we solve Problem 2 implementing the { (weak) greedy} algorithm in the manifold $\lK=\ph^0(\lN)=\{\ph^0_\nu, \nu\in \lN\}$  consisting of minimisers determined by the relation \eqref{phi^0_par}.
The {\it goal} is to  choose $n$ parameters  such that $\Phi_n^0 = {\rm span}\{ \ph_1^0, \ldots, \ph_n^0\}$ approximates the whole manifold $\ph^0(\lN)$ within the given error $\eps$.

To this effect,   as already stated in the introduction, we assume that the matrices   $\mA(\nu)$ and $ \mB(\nu)$ are Lipschitz continuous with respect to the parameter. In turn, this  implies that the mapping  $\nu \to \ph_\nu^0$ possesses the same regularity  as well, with the Lipschitz constant denoted by $C_\ph$.

\subsection{Construction of an approximating space $\Phi_n^0$}

The greedy selection  of each new snapshot relies on the relation \eqref{greedy_error}, which in this setting maximises the distance of elements of $\ph^0(\lN)$ from the space spanned by already chosen  snapshots. Theoretically, this process requires  that we solve \eqref{phi^0_par} for each value of $\nu$.  And this is exactly what we want to avoid. Actually, here we  face the obstacle {\it ii)} from previous section, since  one has to apply the greedy algorithm within a set whose elements are not given explicitly. 
The problem is managed by identifying an appropriate surrogate for the unknown distances  $\dist (\ph^0_\nu, \Phi_n^0)$. To this effect note that
\begin{equation}
  \label{surrogate} 
  \begin{aligned}
\dist (\ph^0_\nu, \Phi_n^0)
&\sim \dist (\Lambda_\nu  \ph^0_\nu, \Lambda_\nu  \Phi_n^0) \cr
&= \dist (x^1 -  e^{T\mA_\nu} \mx^0, \Lambda_\nu  \Phi_n^0) \,,\cr
\end{aligned}
\end{equation}
where $\Lambda_\nu  \Phi_n^0= {\rm span}\{\Lambda_\nu \ph_1^0, \ldots, \Lambda_\nu\ph_n^0\}$, while $\sim$ denotes equivalence of terms resulting from the uniform controllability assumption \eqref{G-bound}.

In such a way we replace the unknown $\ph^0_\nu$ by an easy computed term $x^1 -  e^{T\mA_\nu} \mx^0$, combining the target $x^1$ and the solution of the free dynamics at time $T$. 
As for the other term, note that  $\Lambda_\nu \ph_i^0, i=1...n$ represents  the value at time $T$ of the  solution to the system
\begin{equation*}
\left\{\begin{aligned}
&x'= \mA_\nu x +\mB_\nu u_{\nu,i}\cr
&x(0)=0\cr
\end{aligned}\right.
\end{equation*}
where $u_{\nu,i}$ is the {\it  control} obtained by solving the corresponding adjoint problem (for the parameter $\nu$) with initial datum $\ph_i^0$.

Thus instead of dealing with $\dist (\ph^0_\nu, \Phi_n^0)$ we use the {\it surro\-gate} $\dist (x^1 -  e^{T\mA_\nu} \mx^0, \Lambda_\nu  \Phi_n^0) $, obtained by  projecting an easy obtainable vector $x^1 -  e^{T\mA_\nu} \mx^0$ to a linear space $\Lambda_\nu  \Phi_n^0$ whose basis  is obtained  by solving the adjoint system plus the state one  $n$ times.

\begin{figure}
   \centerline{\includegraphics[scale=0.5]{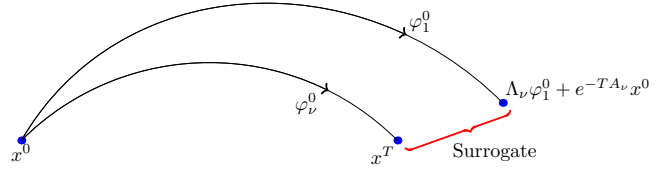}}
 
  \caption{\label{surrogate-graph} The surrogate of $\dist (\ph^0_\nu, \ph_1^0)$ }
\end{figure}

The surrogate measures the control performance of the snapshots $\ph^0_i, i=1..n$ when applied to the system associated to parameter $\nu$  (Figure \ref{surrogate-graph}). Namely,  by relation \eqref{control_nu} the  minimiser $\ph_\nu^0$ uniquely determines the control $u_\nu$ steering the system \eqref{eq1F-d} from the initial datum $x^0$ to the target $x^1$. By replacing $\ph_\nu^0$ with $\ph_i^0$ in \eqref{control_nu} the system is driven to the state $\Lambda_\nu  \ph_i^0 +e^{T\mA_\nu} \mx^0$, whose distance from the target represents the surrogate value.

If  for every $\nu \in \lN$  we can find a suitable linear combination of the above states  close enough  to the target, we   deem that we have found a good approximation of the manifold $\ph^0(\lN)$. Otherwise, we select as the next snapshot  a value for which the already selected snapshots provide  the worst performance. 
The precise description of the offline part of the algorithm is given below.

\vskip 2mm
{\bf Greedy control algorithm - offline part}

\vskip -3mm
\noindent\makebox[\linewidth]{\rule{0.5\textwidth}{0.4pt}}
Fix the approximation error $\eps >0$.

\noindent
    {\bf STEP 1 (Discretisation)}\\    
Choose a finite subset $\tilde \lN$ such that 
$$
\svaki {\nu \in \lN} \quad  \dist(\nu, \tilde \lN) <\delta,
$$
where  $\delta>0$ is a constant to be determined later
  (cf. \eqref{delta-bound}, \eqref{delta}),  in dependence
on the problem under consideration and the tolerance $\eps$. 

\noindent
    {\bf STEP 2 (Choosing $\nu_1$)}\\
   Check the inequality
 \begin{equation}
\label{null_crit}
\max_{\tilde\nu \in \tilde \lN}  |x^1 -  e^{T\mA_{\tilde \nu}} \mx^0|
       < {\eps \over 2} 
\end{equation}
and  stop the algorithm if it holds.
 \\
 Otherwise, determine the first distinguished parameter value as 
    \begin{equation}
      \label{phi_1}
\nu_1
=\argmax_{\tilde \nu \in \tilde \lN} |x^1 -  e^{T\mA_{\tilde \nu}} \mx^0|.
\end{equation}
and choose  $\ph_1^0$ as the minimiser of $J_\nu$ corresponding  to  $\nu=\nu_1$.
\\
\\
\noindent
    {\bf STEP 3 (Choosing $\nu_{j+1}$)}\\
      Having chosen $\nu_1, \ldots, \nu_j$ calculate $\Lambda_{\tilde\nu} \ph_i^0, i=1...j$ for each $\tilde \nu \in \tilde \lN$.
\\
      \noindent
    Check the approximation criteria
 \begin{equation}
\label{stop_crit}
\max_{\tilde\nu \in \tilde \lN}  {\dist (x^1 -  e^{T\mA_{\tilde\nu}} \mx^0, \Lambda_{\tilde\nu}  \Phi_j^0) }
< {\eps \over 2  }\,.
\end{equation}
If the inequality is satisfied, stop the algorithm. 
\\
\noindent
Otherwise, determine the next  distinguished parameter value as
\begin{equation}
\label{step_j}
\nu_{j+1}=\argmax_{\tilde \nu \in \tilde \lN}\dist (x^1 -  e^{T\mA_{\tilde\nu}} \mx^0, \Lambda_{\tilde\nu}  \Phi_j^0)\,,
\end{equation}
choose  $\ph_{j+1}^0$ as the minimiser of $J_\nu$ corresponding  to  $\nu=\nu_{j+1}$ and repeat Step 3.\\
\noindent\makebox[\linewidth]{\rule{0.5\textwidth}{0.4pt}}

The algorithm results in the  approximating space  $\Phi_n^0=\span\{\ph_1^0, \ldots, \ph_n^0\}$, where $n$ is a number of chosen snapshots (specially $\Phi_0^0=\{\vnul\}$ for $n=0$).

 The value of the parameter $\nu_1$ is chosen by testing  the performance of the null control as an initial  guess  for all
  $ \tilde\nu \in \tilde \lN$. 
 The selected value $\nu_1$ is the one for which this performance provides the worst approximation.
 The algorithm is stopped at this initial level only if the null control ensures the uniform control of all system realisations within the given tolerance.

    Note that the set $\{\ph_1^0, \ldots, \ph_n^0\}$ is linearly independent, as  for vectors $\ph_{\tilde\nu}^0$ that linearly depend on already chosen ones, the corresponding  surrogate distance (23), which is the criterion for the choice of new snapshots, vanishes.
Thus the algorithm stops after, at most, $n\leq N$ iterations, and it fulfils the requirements of the weak greedy theory.
More precisely the following result holds.

\begin{theorem}
  \label{main_result}
  Let the $\nu \to (\mA(\nu), \mB(\nu))$ be  Lipschitz and  such that the uniform controllability condition \eqref{G-bound} holds,  and let $C_\ph$ be the Lipschitz constant of  the mapping  $\nu \to \ph_\nu^0$  determined by \eqref{phi^0_par}.

  For a given $\eps>0$ take the discretisation constant $\delta$ such that
 \begin{equation}
    \label{delta-bound}
   \delta\leq \eps/(2C_\ph \Lambda_-).
     \end{equation}
   Then the above algorithm provides a week greedy approximation of the manifold $\ph^0(\lN) $ with the constant
  \begin{equation}
    \label{gamma_const}
    \gamma= \Lambda_- / (2 \Lambda_+) \leq {1 / 2},
  \end{equation}
  and the approximation error less than $\eps/\Lambda_-$. 
  \end{theorem}

\begin{remark}
  The obtained approximation error $\eps/\Lambda_-$ of the family of minimisers $\ph^0(\lN) $ is a consequence of the required $\eps$ approximation of the set $\{\Lambda_\nu \ph_\nu^0, \nu\in \lN\}$ (as $|\Lambda\ph_1^0-\Lambda\ph_2^0|<\eps$ implies $|\ph_1^0-\ph_2^0|<\eps/\Lambda_-$).
    \end{remark}

\begin{remark}
  The case $n=0$,  occurring when inequality \eqref{null_crit} holds,  is a trivial one resulting in a null approximating space that we exclude from the proof. 
  \end{remark}
{\bf Proof:} 
In order to prove the theorem  we have to show:
\begin{itemize}
  \item[{\bf a)}] that  the selected minimisers $\ph_{j}^0$ associated to para\-meter values  determined by \eqref{phi_1} and \eqref{step_j} satisfy
  \begin{equation}
    \label{gamma}
    \begin{aligned}
       |\ph_1^0|
  &\geq\gamma\max_{\nu \in  \lN} |\ph_\nu^0| \cr
       \dist (\ph_{j+1}^0,  \Phi_j^0)
       &\geq \gamma \max_{\nu \in  \lN} \dist (\ph_\nu^0,  \Phi_j^0), \quad j=1..n-1
      \end{aligned}
\end{equation}
  for the constant  $\gamma$ given by \eqref{gamma_const};
\item[{\bf b)}]    that the approximation error 
  $$\sigma_n(\ph^0(\lN)):=\max_\nu\dist(\ph_\nu^0, \Phi_n^0)$$ obtained at the end of the algorithm is less than $\eps/\Lambda_-$. 
\end{itemize}

{\bf a)}  To this effect note that for the first snapshot the following estimates hold:
\begin{equation}
  \label{phi_1-bound}
    \begin{aligned}
  \svaki{ \tilde\nu \in \tilde\lN} \quad |\ph_1^0|
  &\geq{1 \over \Lambda_+} |x^1 -  e^{T\mA_1} \mx^0|\\
&\geq{1 \over \Lambda_+} |x^1 -  e^{T\mA_{\tilde \nu}} \mx^0|
  \geq{\Lambda_- \over \Lambda_+} |\ph_{\tilde \nu}^0|,
  \end{aligned}
\end{equation}

where we have used relation \eqref{phi^0_par} (for $\nu=\nu_1$) and the criterion \eqref{phi_1}.

In order to obtain an estimate including the whole set of parameters $\lN$, we employ  Lipschitz regularity of the mapping  $\nu \to \ph_\nu^0$.
Thus for an arbitrary $\nu \in \lN$, by taking $\tilde\nu\in\tilde\lN$ such that $|\nu -\tilde\nu| < \delta$, by means of \eqref{delta-bound} and \eqref{phi_1-bound} it follows
$$
|\ph_\nu^0|
\leq | \ph_\nu^0 - \ph_{\tilde \nu}^0| + |\ph_{\tilde \nu}^0|
\leq {\eps \over 2 \Lambda_-} + {\Lambda_+ \over \Lambda_-} |\ph_{1}^0|\,.
$$
Having excluded the case \eqref{null_crit}, we have
$$
{\eps \over 2} \leq | \Lambda_1 \ph_1^0| \leq \Lambda_+ |\ph_{1}^0|\,.
$$
implying
\begin{equation}
\label{gamma_1}
|\ph_{1}^0| \geq \gamma \max_{\nu\in\lN} |\ph_\nu^0|\,,
\end{equation}
with $\gamma$ given by \eqref{gamma_const}.

For a general $j$-th iteration,  and an arbitrary $\nu \in \lN$, by taking $\tilde \nu\in B(\nu; \delta)\cap \tilde\lN$ and using \eqref{delta-bound} we get
\begin{equation*}
    \begin{aligned}
      \dist (\ph_\nu^0,  \Phi_j^0)
&\leq | \ph_\nu^0 - \ph_{\tilde \nu}^0| +  \dist (\ph^0_{\tilde\nu},  \Phi_j^0)\\
      &\leq  {\eps \over 2 \Lambda_-}+  {1 \over \Lambda_-}\dist (\Lambda_{\tilde\nu}\ph_{\tilde\nu}^0,  \Lambda_{\tilde\nu}\Phi_j^0)\\
      &\leq  {1 \over \Lambda_-}\left({\eps \over 2}+  \dist (\Lambda_{\tilde\nu}\ph_{j+1}^0,  \Lambda_{\tilde\nu}\Phi_j^0)\right)\,, 
\end{aligned}
 \end{equation*}
where in the last  inequality follows from definition of the next snapshot \eqref{step_j}.  Since the stopping criteria \eqref{stop_crit} is not satisfied until $j<n$, the first term within the last bracket is less than the second one, implying

      \begin{equation*}
    \begin{aligned}
   \dist (\ph_\nu^0,  \Phi_j^0)      &\leq  {2 \over \Lambda_-}\dist (\Lambda_{\tilde\nu}\ph_{j+1}^0,  \Lambda_{\tilde\nu}\Phi_j^0)\\
&\leq 2 {\Lambda_+ \over \Lambda_-} \dist  (\ph^0_{j+1},  \Phi_j^0)\,.
\end{aligned}
 \end{equation*}

 Combining the last  relation with \eqref{gamma_1},  estimate \eqref{gamma} follows.

{\bf b)} Finally, having achieved inequality \eqref{stop_crit} after $n$ iterations,  for an arbitrary $\nu \in \lN$, taking $\tilde\nu$ as above we get
\begin{equation}
  \label{phi_approx}
    \begin{aligned}
\dist (\ph_\nu^0,  \Phi_n^0)
&\leq | \ph_\nu^0 - \ph_{\tilde \nu}^0| +  \dist (\ph_{\tilde\nu}^0,  \Phi_n^0)\\
&< {\eps \over 2\Lambda_-}+{1 \over \Lambda_-}\dist (\Lambda_{\tilde\nu}\ph_{\tilde\nu}^0,  \Lambda_{\tilde\nu}\Phi_n^0)< { \eps\over \Lambda_-}, 
\end{aligned}
 \end{equation}
thus  obtaining the required  approximation error of the set $\ph^0(\lN)$. 
\hfill $\Box$

\begin{remark}

  Besides a choice of  the discretisation constant $\delta$ determined  by \eqref{delta-bound}, other choices are possible as well. Actually, taken any $\delta>0$, define $\tilde C:=  C_\ph \delta/\eps$. The implementation of the algorithm in that case would  lead to the greedy constant $\gamma^{-1}=(2 \tilde C + {1/ \Lambda_-} )\Lambda_+$ and the approximation error of the set $\ph^0(\lN)$ equal to $ \eps (\tilde C + {1/ (2\Lambda_-)})$. Note, however, that no choice of $\delta$ can reduce the order of the error already obtained by \eqref{phi_approx}. \\
 
\end{remark}

\subsection{Construction of an approximate control for a given parameter value }

Having constructed an approximating space $\Phi_n^0$ of dimension $n$, we would like to exploit it for construction of an approximate control $u_\nu^\star$  associated to an arbitrary given value $\nu\in \lN$.

Such a control is given by the relation
\begin{equation}
\label{approx_control}
 u_\nu^\star = \mB_\nu^\ast e^{(T-t) \mA_\nu^*}  \ph_\nu^{0,\star},
\end{equation}

where $\ph_\nu^{0,\star}$ is appropriately chosen approximation of $\ph_\nu^0$  from $\Phi^0_n$.   
It steers the system to the state $x_\nu^\star(T)=\Lambda_\nu \ph_\nu^{0,\star} + e^{T \mA_\nu} x^0$.  Comparing formula \eqref{approx_control} with the one fulfilled by the exact control \eqref{control_nu}, we note that the only difference lies in the replacement of the unknown minimiser $\ph_\nu^0$ by its approximation $\ph_\nu^{0,\star}$.

Thus $x_\nu^\star(T)- e^{T \mA_\nu} x^0$ is a suitable linear combination of vectors $\Lambda_\nu \ph_i^0, i=1..n$. 
As our goal is to steer the system to $x^1$ as close as possible, the best performance is obtained if $\Lambda_\nu \ph_\nu^{0,\star}$ is chosen as projection of $\Lambda_\nu \ph_\nu^0 =x^1- e^{T \mA_\nu} x^0$ to the space $\Lambda_\nu  \Phi_n^0= {\rm span}\{\Lambda_\nu \ph_1^0, \ldots, \Lambda_\nu\ph_n^0\}$.

For this reason we define approximation of $\ph_\nu^0$   as
  \begin{equation}
\label{approx_min}
\ph_\nu^{0,\star}=\sum_i^n \alpha_i \ph_i^0, 
\end{equation}
  where the coefficients $\alpha_i$ are chosen such that $\sum_i^n \alpha_i \Lambda_{\nu} \ph_i^0$ represents  projection of the vector  $x^1 -  e^{T\mA_{\nu}} \mx^0$ to the space $\Lambda_{\nu}  \Phi_n^0$.
  This choice of $\ph_\nu^{0,\star}$ corresponds to the minimisation of the functional $J_\nu$, determined by \eqref{funct_J}, over the 
    space $\Lambda_\nu  \Phi_n^0$.

In order to check performance estimate of the approximate control \eqref{approx_control}, note that for any value $\nu \in \lN$ we have 
\begin{equation}
  \label{perf-est}
     \begin{aligned}
       |x^1-x_{\nu}^\star(T)|&=|\Lambda_{\nu}(\ph_{\nu}^0- \ph_{\nu}^{0,\star})| \cr
  &  \leq  |\Lambda_{\nu}(\ph_{\nu}^0- \ph_{\tilde\nu}^{0,\star})| \cr
       & \leq|\Lambda_{\nu}\ph_{\nu}^0 -  \Lambda_{\tilde\nu}\ph_{\tilde\nu}^0|+ |\Lambda_{\tilde\nu}(\ph_{\tilde\nu}^0- \ph_{\tilde\nu}^{0,\star})|\cr
    &   \quad +  |(\Lambda_{\tilde\nu}- \Lambda_{\nu}) \ph_{\tilde\nu}^{0,\star})|   \,,
\end{aligned}
\end{equation}
where we have used that $\Lambda_{\nu}\ph_{\nu}^{0,\star} $ is the orthogonal projection of $\Lambda_{\nu}\ph_{\nu}^0$ to  the space $\Lambda_{\nu}  \Phi_n^0$ (that also contains $\Lambda_{\nu} \ph_{\tilde\nu}^{0,\star}$), while $\tilde\nu$ is taken from the set $ B(\nu; \delta)\cap\tilde\lN$.

Taking into account the stopping criteria \eqref{stop_crit}, the penultimate term in \eqref{perf-est} is less than $\eps/2$, while the preceding one equals
$$
|(e^{T \mA_\nu} - e^{T \mA_{\tilde\nu}})x^0| \leq C_L |\nu-\tilde\nu| ,
$$
where  $C_L$ denotes the Lipschitz constant of the mapping $\nu \to \Lambda_{\nu}\ph_{\nu}^0=x^1-e^{T \mA_\nu} x^0$.
Similarly, the last term in  \eqref{perf-est} can be estimated from above by $C_\Lambda |\nu-\tilde\nu|$, where $C_\Lambda$ stands for the Lipschitz constant of the mapping $\nu \to \Lambda_{\nu}$.
  Thus in order to  obtain a performance estimate less than   $\eps$, one possibly has  to refine the discretisation used in Theorem \ref{main_result} and to take
\begin{equation}
  \label{delta}
\delta \leq {1 \over 2}  \min\{  {\eps \over C_L+C_\Lambda} ,  {\eps \over  C_\ph \Lambda_-}\},
\end{equation}
where, let it be repeated, $C_\ph$ denotes the Lipschitz constant of the mapping $\nu \to \ph_\nu^0$.

This finalises solving of Problem 2 and leads us to the following result. 

\begin{theorem}
  \label{thm-online}
  Let $\nu \to (\mA(\nu), \mB(\nu))$ be a Lipschitz mapping such that the uniform controllability condition \eqref{G-bound} holds.

 Given $\eps>0$,  let  $\Phi_n^0 = {\rm span}\{ \ph_1^0, \ldots, \ph_n^0\}$ be an approximating space constructed by the greedy control algorithm with the discretisation constant $\delta$ given by \eqref{delta}. Then for any $\nu \in \lN$ the approximate control $u_\nu^\star$ given by \eqref{approx_control}  and \eqref{approx_min} steers  the control system \eqref{eq1F-d} to the state $x_\nu^\star(T)$ within the $\eps$ distance from the target $x^1$. 
 \end{theorem}
\begin{remark}
  As already commented in Section \ref{control-prel}, the approximate control $u_\nu^\star$ does not belong to the space $\span\{u_1,\ldots, u_n\}$ spanned by controls associated to selected parameter values, as the control operator and system matrix entering \eqref{approx_control} correspond to the given value $\nu$, while $u_i = \mB_i^\ast e^{(T-t) \mA_i^*}  \ph_i$.
\end{remark}

\begin{remark}
  The greedy control also applies if we additionally assume that initial datum $x^0$, as well as target state $x^1$ depend on the parameter $\nu$ in a Lipschitz manner. In that case, the greedy search is performed in the same manner as above, i.e.  by exploring   elements of the set $\{x^1_\nu  -  e^{T A_\nu} x^0_\nu |\nu \in \lN\}$ in calculation of the  surrogate distance, and the obtained results remain valid with the same constants.
\end{remark}

\begin{remark}
  \label{cont_dep}
    All the above results on greedy control remain valid if  instead of Lipschitz continuity we  merely assume continuous dependence with respect to the parameter $\nu \in \lK$. Namely, as $\lK$ is a compact set, the assumption directly implies uniform continuity, which suffices for the proof of Theorems \ref{main_result} and \ref{thm-online} in which we need $\ph_\nu^0 $ and $\ph_{\tilde\nu}^0 $ to be close whenever $\nu$ and $\tilde\nu$ are. The only difference in that case is that  the discretisation constant $\delta$ can not be given explicitly in terms of $\eps$, unlike expressions \eqref{delta-bound} and \eqref{delta}.
\end{remark}

We finish the section describing the algorithm  summarising the above procedure to construct the approximate control $u_\nu^\star$.
Hereby we suppose that, for fixed $\eps>0$, the approximating space   $\Phi_n^0 = {\rm span}\{ \ph_1^0, \ldots, \ph_n^0\}$ has been constructed   by the offline part of the greedy control algorithm  with the  discretisation  constant $\delta$ given by \eqref{delta}.

\vskip 2mm
{\bf Greedy control algorithm - online part}

\vskip -3mm
\noindent\makebox[\linewidth]{\rule{0.5\textwidth}{0.4pt}}
A parameter value $\nu \in \lN$ is given.

\noindent
    {\bf STEP 1 }  
    Calculate $\Lambda_\nu \ph_\nu^0 =x^1- e^{T \mA_\nu} x^0$.

    {\bf STEP 2 }
    Calculate $\Lambda_\nu \ph_i^0=e^{(T-t) \mA_\nu} B u_i^\star,\; i=1..n$,
      where $
      u_i^\star=\mB_\nu^\ast e^{(T-t) \mA_\nu^*}  \ph_i^{0}.
     $
    
    {\bf STEP 3 }
    Project $\Lambda_\nu \ph_\nu^0$ to $\Lambda_\nu  \Phi_n^0= {\rm span}\{\Lambda_\nu \ph_1^0, \ldots, \Lambda_\nu\ph_n^0\}$.
    Denote the projection by $P_n \Lambda_\nu \ph_\nu^0$.

    {\bf STEP 4}
    Solve the system
    $P_n \Lambda_\nu \ph_\nu^0= \alpha_i \Lambda_{\nu} \ph_i^0$ for $\alpha_i, i=1..n$.

    {\bf STEP 5}
    The approximate control is given by
    $$
    u_\nu^\star=\sum_i \alpha_i u_i^\star,
    $$
    where $u_i^\star$ are already determined within Step 2. 
    \vskip -3mm
    \noindent\makebox[\linewidth]{\rule{0.5\textwidth}{0.4pt}}

For any $\nu \in \lN$ the obtained approximate control steers the system \eqref{eq1F-d} within an $\eps$ distance from the target.

Note that for a parameter value $\tilde\nu$ that  belongs to the discretisation set $\tilde \lN$ used in the construction of an approximating space $\Phi_n^0$
steps 1--3 of the last algorithm can be skipped as the corresponding terms have already been calculated within the offline part of the algorithm.

\section{Computational cost}
\subsection{Offline part}

The offline part of the greedy control algorithm consists of two main ingredients. First, the search for distinguished parameter values $\nu_j$ by examining the surrogate value over the set $\tilde N$, and, second,  the calculation of the corresponding  snapshots $\ph_j^0$.
This subsection is devoted to estimate its computational cost. 

\noindent
    {\bf Choosing $\nu_1$}\\  
In order to identify the first  parameter value,  one has to maximise the expression $|x^1 -  e^{T\mA_{\tilde \nu}} \mx^0|$ over the set $\tilde\lN$, which represents the distance of the target $x^1$ from the free dynamics state. To this effect one has to solve  the original system \eqref{eq1F-d} with zero control  $k=\card(\tilde\lN)$ times. Denoting the cost of solving the system by $C$, the cost of choosing $\nu_1$ thus equals
$$
k (C + 3N),
$$
where the second term corresponds to calculation of the distance between  two vectors in $\R^N$.

In general, the cost $C$ differs depending on the type of matrices $\mA_\nu$ under consideration and the method chosen for solving the corresponding control system. For example,  an implicit one step method consists on solving $T/\Delta t$ linear systems, where $\Delta t$ denotes the time discretisation step. However, as we consider time independent  dynamics, all these  systems have the same matrix, thus LU factorisation is required just ones.
Consequently, the cost $C$ can be estimated as
\begin{equation}
  \label{nu_1}
C={2 \over 3} N^3 + 
{T \over \Delta t} 4 N^2\,,
\end{equation}
where the first part corresponds to the $LU$ factorisation, while the latter one is the cost of building and solving a system of a type $LU x_{k+1}= b_k$.

\noindent
    {\bf Calculating $\ph_1^0$}\\
In order to determine the first snapshot one has to solve the system \eqref{phi^0_par} for the chosen parameter value $\nu_1$. To this effect we  construct (the unknown) Gramian matrix $\Lambda_1$ by calculating $\Lambda_1 e_i, , i=1...N$ for   the vectors of the canonical basis $e_i$. According to the results presented in Section \ref{control-prel}, for any  $\ph^0 \in \R^N$ the corresponding vector $\Lambda_1 \ph^0$ can be determined as the state of the control system at time $T$ with the control $u=\mB^\ast \ph$, with $\ph$ being the solution of the adjoint problem starting from $\ph^0$.

In such a way, the cost of  composing the system for $\ph_1^0$ equals $2N C$ (corresponding to  $N$ control plus adjoint problems starting from $N$ different data).  As solving the system requires  a number of operations of order $2/3 N^3$, the cost of this part of the algorithm turns out as
$$
2N C + {2 \over 3} N^3\,.
$$

\noindent
    {\bf Choosing $\nu_{j+1}$}\\
Suppose we have determined the first $j$ snapshots and have constructed the approximating space $\Phi_j^0=\span \{ \ph_1^0, \ldots, \ph_j^0\}$. The next parameter value is chosen by  maximising the distance of the vector $x^1 -  e^{T\mA_{\tilde\nu}} x^0$, calculated already in the first iteration, from the space $\Lambda_{\tilde\nu}  \Phi_j^0$ over the set $\tilde\lN\setminus \{ \nu_1,  \ldots , \nu_j\}$ (the already chosen values can obviously be excluded from the search).  The basis of the above space has been determined gradually throughout the previous iterations up to the last vector $\Lambda_{\tilde\nu}  \ph_j^0$, whose calculation is performed at this level. As explained above, it  requires  $2 C$ operations, that have to be performed $k-j$ times.

In addition, this basis is orthonormalised, thus enabling the efficient calculation of  the distance in \eqref{stop_crit}. This process is performed gradually throughout the algorithm as well, and in each iteration we just  orthonormalise the last vector with respect to the the rest of the, already orthonormalised set. The corresponding cost for a single value $\tilde\nu \in\tilde\lN$ is of order
$4N j$. Similarly, the (orthogonal) projection of
 $x^1 -  e^{T\mA_{\tilde\nu}} x^0$ to $\Lambda_{\tilde\nu}  \Phi_j^0$ takes into account its projection to $\Lambda_{\tilde\nu}  \Phi_{j-1}^0$ used in the previous iteration, and adds just a projection to the last introduced vector.

As its corresponding cost is of order $4N$, the total cost of this part of the algorithm equals
$$
(k-j)(2 C + 4N j+4N +3N)\,,
$$
where, as in \eqref{nu_1}, the last term corresponds to calculating the distance between  two vectors in $\R^N$. 

\noindent
    {\bf Calculating $\ph_{j+1}^0$}\\
    In order to determine the next snapshot, we have to construct corresponding Gramian matrix $\Lambda_{j+1}$. As explained in the part related to the calculation of $\ph_1^0$, this can be done by applying the Gramian to  some basis of $\R^N$. As in the previous part of the algorithm we have already calculated $\Lambda_{j+1} \ph_i^0, i=1...j$, it is enough to calculate $\Lambda_{j+1} \e_{j_i}, i=1...N-j$, for vectors of canonical basis $\e_{j_i}$ complementing the set $\Phi_j^0$ to the basis of $\R^N$.   Thus building the matrix of the system \eqref{phi^0_par} requires $2(N-j)C$ operations.

    In addition, the complementation of the set $\Phi_j^0$ takes advantage of the basis obtained as complementation of the set $\Phi_{j-1}^0$ in the previous iteration.
  Adding a new snapshot $\ph_j^0$ to that basis,  one of  vectors $e_{{j}_i}$ has to be removed such that the new set results in a basis again.
  Identifying this vector requires solving  a single $N\times N$ system.

   Finally, taking into account the cost of solving the system \eqref{phi^0_par}, the cost of calculating the next snapshot turns out to be
$$
2(N-j)C + {4 \over 3} N^3\,.
$$

\noindent
    {\bf Total cost}\\
Summing the above costs for $j=1...n$ the total cost of the algorithm results in
\begin{equation}
  \label{total}
  \begin{aligned}
k(&C+3N) + n(k-{n\over 2})(2C+7N)+2n^2N(k-{2 \over 3}n) \\
&+ 2C n (N-{n \over 2}) + {4 \over 3} N^3 n
\end{aligned}
\end{equation}

As the cost $C$ of solving the control system is $O(N^3)$, the most expensive part of the greedy control algorithm corresponds to the terms containing this cost.

It is interesting to notice that, as the number of chosen snapshots $n$ approaches either the number of eligible parameter values $k=|\tilde\lN|$, or the system dimension $N$, this part converges to
$$
C(k+2k N),
$$
which corresponds, up to lower order terms, to the cost of calculating $\ph_{\tilde \nu}^0$ for all $k$  values  $\tilde\nu \in \tilde\lN$.

This demonstrates that the application of the greedy control algorithm is always cheaper than a naive approach that consists of calculating controls for all values of the parameter from  an uniform mesh on $\lN$, taken fine enough so to achieve the approximative control \eqref{eq2F-dap}.

Let us note that in the above analysis we have assumed that in the each iteration  we calculate the surrogate distance \eqref{stop_crit} for all values of the set $\tilde\lN$, except those already chosen by the greedy algorithm. And that we  repeat the procedure until the residual is less than $\eps/2$ for all parameter values of  $\tilde\lN$.
However,  very likely, for some, or even many parameter values, the corresponding residual will be less than the given tolerance already before the last iteration $n$. And  once it occurs,   these values do not have to be explored any more, neither  we  have to calculate the corresponding surrogates for subsequently chosen snapshots.
Therefore,  the  obtained cost estimates are rather conservative, and in practice we expect the real cost to be lower than the above one.

Further significant cost reductions are obtained under assumption of  a system  with  parameter independent matrix $\mA$ and the control operator $\mB$ (with dependence remaining only in initial datum and/or in the target). In that case it turns that  the corresponding Gramian matrix $\L$ is also parameter independent.  For this reason,  maximising the distance appearing in \eqref{stop_crit} requires solving of control plus adjoint problem just once in each iteration. As a result, the most expensive term of the total cost \eqref{total}  is replaced by quite a  moderate one
$$
C(k+2n+2N)\,.
$$

\subsection{Online part}
\noindent
    {\bf STEP 1}\\
Calculating $\Lambda_\nu \ph_\nu^0$ corresponds to solving the control system once, whose computational cost equals $C$. 

\noindent
    {\bf STEP 2}\\
    Calculating $\Lambda_\nu \ph_i^0$ requires solving the loop of the adjoint plus control system $n$ times, resulting in the cost of $2nC$.

    \noindent
    {\bf STEP 3}\\
    The most expensive part of this step consist in the  Gram Schmidt orthonormalisation procedure of the set  $\{\Lambda_\nu \ph_1^0, \ldots, \Lambda_\nu \ph_n^0\}$, whose cost equals $2Nn^2$.

    \noindent
    {\bf STEP 4}\\
    Solving the system with QR decomposition requires $2N n^2 - 2 n^2/3$ operations.

   \noindent  
 {\bf STEP 5}\\
    The cost of this step is negligible compared to the previous ones. 

  \noindent  
    {\bf Total cost}\\
    Thus we obtain that  the total cost of finding an approximative control $u_\nu^\star$ for a given parameter value equals
    $$
    C(1+2n) + 4 N n^2 - 2 n^2/3.
    $$
  As $n$ approaches the system dimension  $N$, its most expensive part  $C(1+2n)$   converges, up to lower order terms,  to the cost of calculating an exact control $u_\nu$.   Consequently, the cost reduction obtained by choosing the approximative control $u_\nu^\star$  obtained by the greedy control algorithm instead of the exact one depends linearly on the ratio between the number of used snapshots $n$ and the system dimension $N$.

\section{Infinite dimensional problems}
\label{infinite}

The theory and the (weak) greedy control  algorithm developed in the preceding section for finite dimensional linear control systems extend to ODEs in infinite dimensional spaces. They can be written exactly as in the form (\ref{eq1F-d}) except for the fact that solutions $x(t, \nu)$ for each value of the parameter $\nu$ and each time $t$ live in an infinite dimensional Hilbert space $X$. The key assumption that distinguishes these infinite-dimensional ODEs from Partial Differential Equations (PDE)  is that the operators $(\mA(\nu), \mB(\nu))$ entering in the system are assumed to be bounded.

The controllability of such systems has been extensively elaborated during last few decades, expressed in terms of semigroups generated by a bounded linear operator, see, for instance,\cite{Tri,TuW,Zab}. In fact, the existing theory of controllability for linear Partial Differential Equations (PDE) can be applied in that context too. This allows to characterise controllability in terms of the observability of the corresponding adjoint systems. In this way, the uniform controllability conditions of parameter dependent control problems can be recast in terms of the uniform observability of the corresponding adjoint systems.

However, in here, we limit our analysis to the case of infinite-dimensional ODEs for which the evolution is generated by bounded linear operators, contrarily to the case of PDEs in which the generator is systematically an unbounded operator. The greedy theory developed here is easier to implement in the context of infinite-dimensional ODEs  since, in particular,  the analytic dependence property of controls with respect to the parameters entering in the system can be more easily established.

In fact, in the context of  infinite-dimensional ODEs, most of the results in Section 2 on finite-dimensional systems apply as well. In particular, the characterisation of the controls as in (\ref{control_nu}), in terms of minimisers of quadratic functionals of the form $J$ as in (\ref{funct_J}) holds in this case too, together with (\ref{Gramian}) for the Gramian operators. Obviously, the Kalman rank condition cannot be extended to the infinite-dimensional case. But the open character of the property of controllability with respect to parametric perturbations remains true in the infinite-dimensional case too, i.e. if the system is controllable for a given value of $\nu$, under the assumption that $(\mA(\nu), \mB(\nu))$ depends on $\nu$ in a continuous manner, it is also controllable for neighbouring values of $\nu$.

In this section we shall analyse convergence rates of the constructed greedy control algorithm as the dimension of the approximating space tends to infinity, an issue that only makes sense for infinite-dimensional systems. In particular, we are interested in the dimension of the approximating space $\ph^0(\lN)$ required to provide an uniform control within the given tolerance $\eps$.  This problem can be stated in terms of the estimate \eqref{stop_crit}:  what is the number of the algorithm iterations $n$ we have to repeat until the estimate is fulfilled. 

In case of systems  of a finite dimension $N$, the algorithm constructs an $n$ dimensional approximating space of $\ph^0(\lN) \subseteq \R^N$, and it  certainly stops after, at most, $n\leq N$ iterations. For infinite-dimensional systems, however there is no such an obvious stopping criteria.  

In general we  analyse the performance of the algorithm by comparing  the greedy approximation errors $\sigma_n(\ph^0(\lN))=\max_\nu\dist(\ph_\nu^0, \Phi_n^0)$ with the Kolmogorov widths $d_n(\ph^0(\lN))$, which represent the best possible approximation of $\ph^0(\lN)$ by a $n$ dimensional subspace in $\R^N$.
To this effect one could try to employ Theorem \ref{greedy_rates} which connects sequences $(\sigma_n)$ and $(d_n)$. However,  we have to apply the theorem to the set  $\ph^0(\lN)$ (instead of $\lN$). But while  Kolmogorov width of a set of admissible parameters $\lN$ is usually easy to estimate, that is not the case for a corresponding set of  solutions (to a parametric dependent equation), or  minimisers (as studied in this manuscript).  Fortunately, a result in that direction has been provided recently for holomorphic mappings (\cite{CDV}) under the assumption of a polynomial decay of Kolmogorov widths. 
\begin{theorem}
  \label{greedy_rates2}
  For a pair of complex Banach spaces $X$ and $V$ assume that $u$ is a holomorphic map from an open set $\lO\subset X$ into $V$ with uniform bound. If $\lK\subset\lO$ is a compact subset of $X$ then for any $\alpha>1$ and $C_0>0$
  \begin{equation}
  \label{holo}
    d_n(\lK)\leq C_0 n^{-\alpha} \quad \Longrightarrow \quad d_n(u(\lK)) \leq C_1 n^{-\beta}, \quad n\in\N,
    \end{equation}
    for any $\beta<\alpha-1$ and the constant $C_1$ depending on $C_0$, $\alpha$ and the mapping $u$. 
\end{theorem}
  \begin{remark}
  The proof of the theorem provides an explicit estimate of the constant $C_1$ in dependence on  $C_0$, $\alpha$ and the mapping $u$.  However, due to its rather complicated form we do not expose  it here.
  \end{remark}
  
Going back to our problem, Theorem \ref{greedy_rates2} can be applied  under the assumption that the  mapping $\nu \to \big(\mA(\nu), \mB(\nu)\big)$ is analytic (its image being embedded in the space of linear and bounded operators in $X$), which implies, in view of the representation formula (\ref{control_nu}),  that the  mapping  $\lN$ to $\ph^0(\lN)$ is analytic as well. Note that this issue is much more delicate in the PDE setting, with the generator of the semigroup $\mA(\nu)$ being an unbounded operator. In fact the property fails in the context of hyperbolic problems although it is essentially true for elliptic and parabolic equations.

As we consider a finite number of  parameters, lying in the set $\lN\subset\R^d$, the polynomial decay of $d_n(\lN)$ can be achieved at any rate $\alpha>0$ just by adjusting the corresponding constant $C_0=C_0(\alpha)$ in \eqref{holo} (note that $d_n(\lN)=0$ for $n\geq d$). Of course, the Kolmogorov widths of the set $\ph^0(\lN)$ do not have to vanish for $n$ large, but Theorem \ref{greedy_rates2} ensures their polynomial decay at any rate.

 Combining Theorems \ref{greedy_rates}, \ref{main_result} and  \ref{greedy_rates2} we  thus obtain the following result.
 \begin{corollary}
   \label{cor-decay}
    Let  the mapping $\nu \to (\mA(\nu), \mB(\nu))$, corresponding to the parameter dependent infinite-dimensional ODEs, be  analytic and such that the uniform controllability condition \eqref{G-bound} holds.
Then the greedy control algorithm  ensures a polynomial decay of arbitrary order of the approximation rates.

More precisely,  for all $\alpha >0$ there exists $C_\alpha >0$ such that for any $\nu$ the  minimiser $\ph^0_\nu$ determined by the relation \eqref{phi^0_par} can be approximated by linear combinations of the weak-greedy ones as follows:
$$
dist( \ph^0_\nu; \span\{\ph^0_j: j=1,..., n\}) \le C_\alpha n^{-\alpha},
$$
where  $C_\alpha$ can be determined by exploring constants appearing in \eqref{poly},  \eqref{gamma_const}  and  \eqref{holo}.
 \end{corollary}
 
   The last result provides us with a stopping criteria of the greedy control algorithm. For a given tolerance $\eps>0$, the algorithm  stops  after choosing $n>\root \alpha \of{C_\alpha/\eps}$ snapshots. It results with a $n$ dimensional space  $\Phi_n^0$  approximating the family of minimisers $\ph^0(\lN)$ within the error $\eps/\Lambda_-$, and providing,  by means of formul\ae \, \eqref{approx_control}  and \eqref{approx_min},  a uniform control of our system within an  $\eps$ tolerance.

A similar result holds if $\lN$ is infinite-dimensional, provided  its Kolmogorov width decays polynomially. A typical example of such a set is represented by the so called affine model in which the parameter dependence is given by
\begin{equation}
\label{affine}
   \begin{aligned}
\mA( \mnu)= \sum_{l=1}^\infty \nu_l \mA_l,
\end{aligned}
\end{equation}
and/or similarly for $\mB(\mx, \mnu)$.
Here it is assumed that $\mnu$ belongs to the unit cube in $l^\infty$, i.e. that $|\nu_l| \leq 1$ for any $l\in \N$, while the sequence of numbers $a_l:=||\mA_l||$  belongs to $l^p$ for some $p<1$.

However, note that in this case the Kolmogorov width of the set $\lN =B(0;1, l^\infty)$ does not decay polynomially (actually these are constants equal to 1), but the polynomial decay is obtained for the set $\mA(\mnu)= \{\mA(\cdot, \mnu); \mnu \in \lN\}$. Indeed, rearranging the indices so that  the sequence $(a_l)$ is decreasing, it follows
\begin{equation*}
  \begin{aligned}
||\mA( \mnu)- \sum_{l=1}^n \nu_l \mA_l||
&\leq \sum_{l=n+1}^\infty \nu_l a_l \\
&\leq a_{n+1}^{1-p} \sum_{l=n+1}^\infty \nu_l a_l^p
\leq C n^{1-1/p} \,
\end{aligned}
\end{equation*}
where we have used that for a decreasing $l^p$ sequence we must have $a_n^p \leq C n^{-1}, n\geq 1$.

Thus one can consider the mapping $\big(\mA(\mnu), \mB(\mnu)\big) \to \ph_\mnu^0$ and apply Theorem \ref{greedy_rates2} to $K=\big(\mA(\mnu), \mB(\mnu)\big)$.

Furthermore, in the case of the affine model this theorem can be improved, implying the  Kolmogorov $n$-widths of sets $K$ and $u(K)$ decay at the same rate $n^{-\alpha},  \alpha= 1/p -1$ (e.g. \cite{CDV}) . 
Consequently, one obtains that the greedy approximation rates $\sigma_n(\ph^0(\lN)$ decay at the same rate as well.

  Finally, let us mention that the  cost of the greedy control is significantly reduced if one considers an affine model in which the  control operator has finite  representation of the form \eqref{affine}, while the system matrix $\mA$ is taken as  parameter independent.

In that case the corresponding Gramian is of the form:
   \begin{equation}
      \label{affine_Lambda}
  \Lambda(\mnu)=\sum_{l,m =1}^L \nu_l\nu_m \Lambda_{lm} \,
  \end{equation}
   where $\Lambda_{lm}=\int_0^T e^{(T-t)\mA}\mB_l \mB_m^\ast e^{(T-t)\mA^\ast} dt$, while $\mnu=(\nu_1, \ldots, \nu_L)$.
Here we consider  a finite representation, but a more general one can be reduced   to this one by truncation of the series. 
Consequently,  computing $\Lambda_{\tilde\mnu}  \ph_j^0$ for a chosen snapshot does not require solving the loop of the adjoint plus control system for each value $\tilde\mnu \in \tilde\lN$. Instead, it is enough to perform the loop $L^2$ times in order to obtain vectors $\Lambda_{lm}  \ph_j^0, l,m=1..L$, and express $\Lambda_{\tilde\mnu}  \ph_j^0$ as their linear combination by means of \eqref{affine_Lambda}.

Such approach can  result in a lower  cost of the greedy control algorithm compared to the one obtained in the previous section,
with a  precise  reduction rate depending on the relation between the series length $L$ and the number of eligible parameters $k$. 

\section{Numerical examples}

\subsection{Wave equation}

We consider the control system \eqref{eq1F-d} whose governing matrix has the  block form 
  $$
  \mA=\left(\begin{matrix}
  \mnul & \mI_{N/2}\cr
  \nu ( N/2+1)^2 \tilde \mA & \mnul\cr
  \end{matrix}\right),
  $$
  where $\mI_{N/2}$ is the identity matrix of dimension $N/2$, while 
  \begin{equation}
      \label{A}
    \tilde \mA=\left(\begin{matrix}
    -2 & 1& 0 & \cdots &0& 0\cr
    1&  -2 & 1&  \cdots &0& 0\cr
    0 & 1&  -2 &  \cdots &0& 0\cr
    \vdots&\vdots&\vdots&\ddots&\vdots&\vdots\cr
    0 &0&  0 & \cdots& -2 &1\cr
    0 &0&  0 & \cdots& 1  &-2\cr
    \end{matrix}\right).
    \end{equation}
    The control operator is assumed to be of the parameter-independent form
    $$
    \mB= 
    (0, \ldots,0, ( N/2+1)^2 )^\top.
   $$

    The system corresponds to the semi-discretisation of the wave equation problem with the control on the right boundary:
    \begin{equation}
      \label{wave}
      \left\{         \begin{split}
    \partial_{tt} v -  \nu \partial_{xx} v& =0, \qquad (t, x) \in \oi 0T \times \oi 01,\\
     v(t, 0)& =0, \quad v(t, 1)=u_\nu(t),\\
      v(0, x)& =v_0(x), \quad \partial_t v(x, 0)=v_1(x)\,.\\
     \end{split}\right.
  \end{equation}
     Parameter $\nu$ represents the (square of ) velocity of propagation, while     $N/2$ corresponds to the number of inner points in the discretisation of the space domain.

   For this example we specify the following ingredients:  
  $$T=3,\;  N=50, 
  v_0=\sin(\pi x),  v_1=0 \,.
  $$
  The final target is set at $
  x^1=\vnul,  $
    and we assume
  $$
  \nu \in \zi 1{10} = \lN .
  $$
  The system satisfies the Kalman's rank condition and accordingly the control exists for any value of  $\nu$. Although the convergence of this direct approximation method, in which one computes the control for a finite-difference discretisation, hoping that it will lead to a good approximation of the continuous control, is false in general (see \cite{Z1}), this is a natural way to proceed and it is interesting to test the efficiency of the greedy approximation for a given $N$, which in practice corresponds to fixing the space-mesh for the discretisation. In any case, the time-horizon $T$ is chosen such that the Geometrical Control Condition (ensuring controllability of the continuous problem \eqref{wave}, see \cite{BLR}) is satisfied  for the continuous wave equation and all  $\nu\in \lN$. 
  
  The greedy control algorithm has been applied with $\eps=0.5$ and the uniform discretisation of $\lN$ in $k=100$ values.
  
\begin{figure}
          \vskip -4cm  
     \centerline{  \includegraphics[scale=0.5]{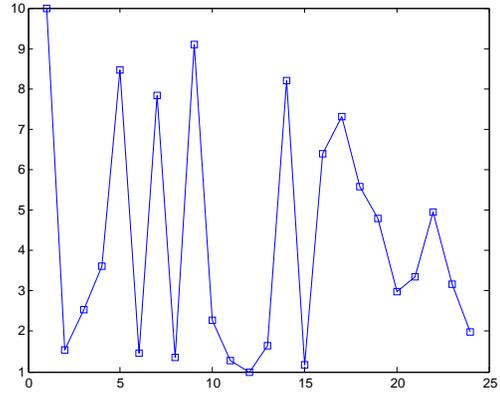}}
      \vskip -4cm
  \caption{\label{par_values} Distribution of the  selected parameter values. }
  \end{figure}

  The algorithm stopped after 24 iterations, choosing 
  24 (out of 100) parameter values
  depicted on Figure \ref{par_values}.
    The figure manifests the way by which  the algorithm selects the values, taking them from different parts of the parameter set in a zigzag manner, in order to obtain the best possible approximation performance. 
      The corresponding minimisers $\ph_i^0$ have been calculated and stored, completing the offline part of the process.
     The elapsed time of the algorithm is 2\,312 seconds, performed on a personal notebook with a 2.7 GHz processor and DDR3 RAM with 8 GB and 1,6 GHz.

   In the online part one explores the stored data in  order to construct approximate controls  for different parameter values ranging between 1 and 10  by means of formula \eqref{approx_control}. Here we present the results obtained for $\nu=\pi$.
   
\begin{figure}
      \vskip -4cm  
     \centerline{  \includegraphics[scale=0.5]{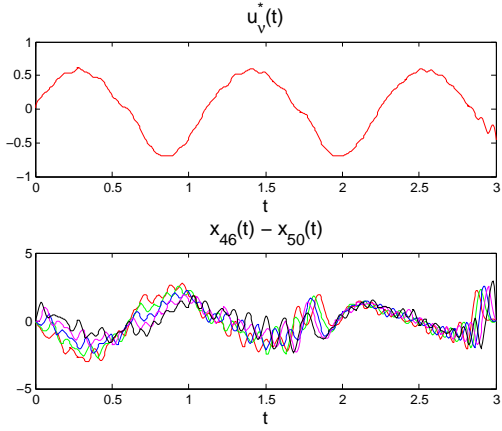}}
      \vskip -4cm
  \caption{\label{wave-control}Evolution of a) the approximate control and 
     b) the last 5 system components   for 
     $\nu=\pi$.}
  \end{figure}

   The approximate control constructed by the greedy control algorithm  is depicted on  Figure \ref{wave-control}.a. 
It is worth to mention that the constructed control turns out to  be  almost indistinguishable from the optimal exact control obtained by standard methods, based on minimisation of the  cost functional  $J$ from \eqref{funct_J}. Although providing almost the same output as the classical methods, the advantage of the approach presented in this paper lies in the gain on computational cost.
 The elapsed time of the online part of the algorithm in this case is 7 seconds, compared to 51 second required to construct the optimal exact control by standard methods. This result confirms the computational efficiency of the greedy approach elaborated in Section 5.

    Figure \ref{wave-control}.b displays  the evolution of the last 5 components of the system, corresponding to the time derivative of the solution to \eqref{wave} at grid points, controlled by the approximate control $u_\nu^\star$.  As required, their trajectories are driven (close) to  the zero target.

\begin{figure}
      \vskip -4cm  
      \centerline{ \includegraphics[scale=0.5]{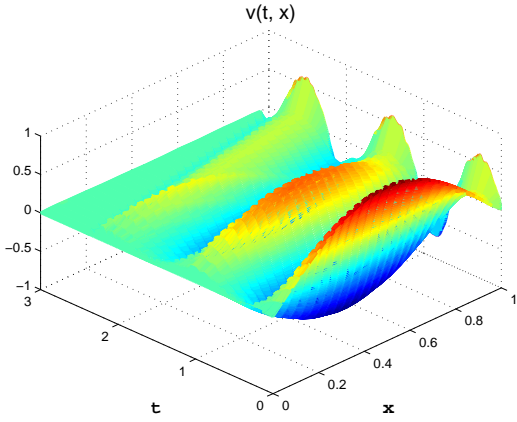}}
      \vskip -4cm  
        \caption{\label{wave-3D} Evolution of the solution to the semi-discretised problem  \eqref{wave} governed by the approximate control $u_\nu^\star$ for $\nu=\pi$. }
        \end{figure}

   The first 25 components of the system are plotted by a 3-D plot (Figure \ref{wave-3D}) depicting the evolution of the solution to a semi-discretised problem \eqref{wave} governed by the approximate control $u_\nu^\star$.
       The starting wave front, corresponding to the sinusoidal initial position, as well as the oscillations in time, gradually diminish and the solution is steered (close) to zero, as required. Furthermore, the total distance of the solution from the target equals $|x(T)|=0.05$.

\begin{figure}
     \vskip -1cm  
      \centerline{   \includegraphics[scale=0.25]{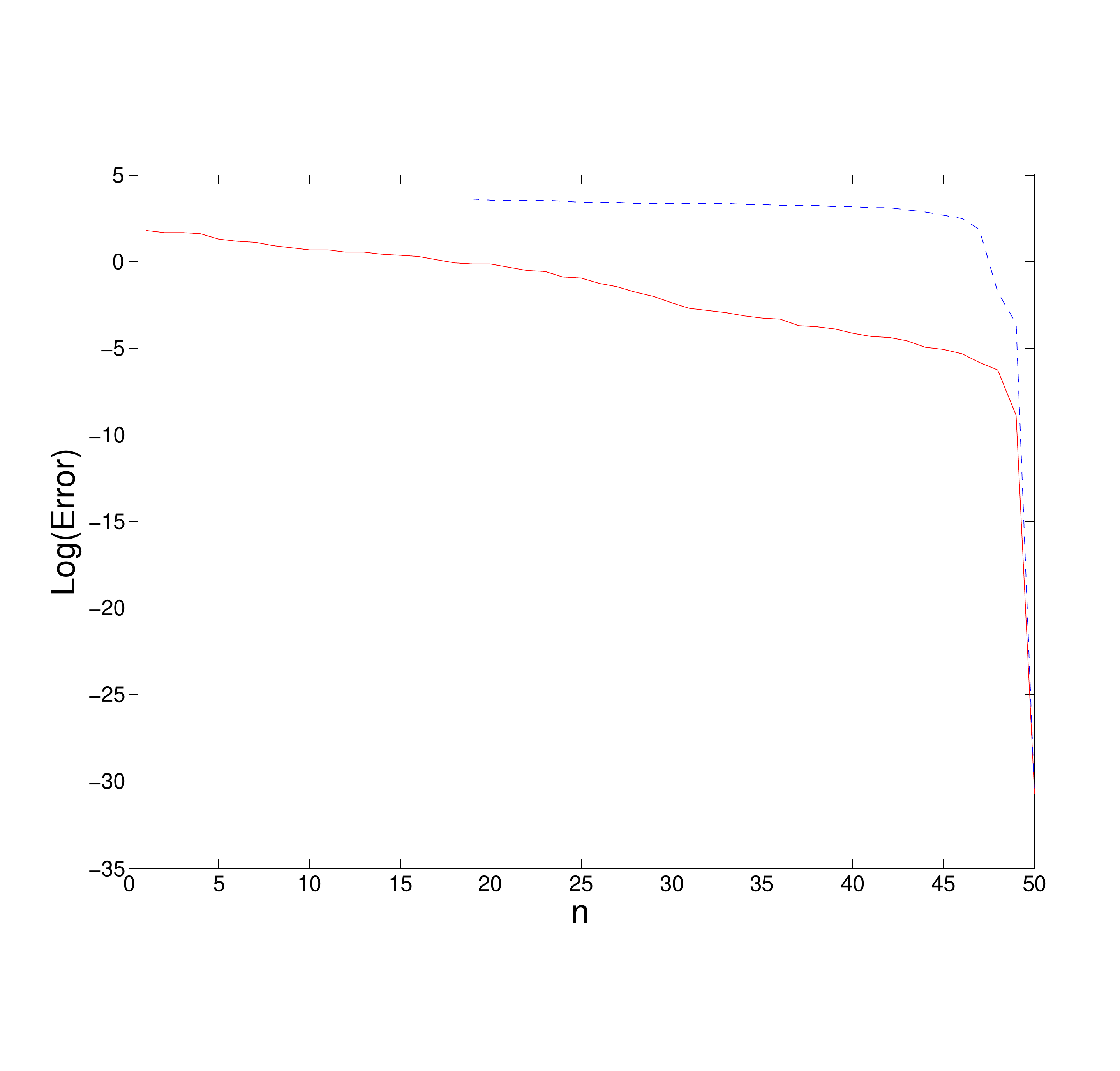}}
    \vskip -1cm  
        \caption{\label{approx_rates}Greedy (solid) vs canonical (dashed) approximation rates.}
  \end{figure}

        The efficiency of the greedy control approach is checked by exploring the corresponding approximation rates $\sigma_n(\ph^0(\tilde\lN))=\max_{\tilde\nu}\dist(\ph_{\tilde\nu}^0, \Phi_n^0)$,
                  whose logarithmic values are depicted by the solid curve of Figure \ref{approx_rates}. The curve decreases and goes below $\log(\eps)=\log(0.5)$ for $n=24$,  stopping the algorithm after 24 iterations.
          
        Almost linear shape of the curve suggests an exponential decay of the 
approximation rates $\sigma_n(\ph^0(\tilde\lN))$. Such appearance is in accordance with Corollary \ref{cor-decay} providing their polynomial decay at any rate. 
   The abrupt breakup of the linear behaviour near $n=50$ is due to the fact that we consider a finite dimensional problem in which the approximation space is saturated when  $n$ reaches the size of the ambient space $\R^N$.    
        
         The bases of the approximation spaces $\Phi_n^0$ are built iteratively in a {\it greedy} manner, by exploring the surrogate distance \eqref{surrogate}  over the set $\tilde N$. As explained in the previous section, this search is costly, but should produce optimal approximation rates in the sense of the Kolmogorov widths. As opposed to that, one can explore some a cheaper approach, in which vectors spanning  an approximation space are  chosen arbitrarily, e.g. by taking vectors of the canonical basis.

         The dashed curve of Figure \ref{approx_rates} plots approximation errors of   spaces $E_n$ spanned by the first $n$ vectors of the canonical basis in $\R^N$. Obviously the greedy approach wins over the latter one, in accordance with the theoretical results, ensuring  optimal performance.
         Eventually, both curves exhibit sharp decay at the moment in which the size of the approximation space $n$ reaches the size of the ambient space $\R^N$. 

\subsection{Heat equation}

For  $   \mA=( N+1)^2\tilde\mA    $, with $\tilde\mA$ given by  \eqref{A}, and the control operator $
    \mB= 
    (0, \ldots,0, ( N+1)^2 )^\top,
   $
  the system \eqref{eq1F-d} corresponds to the space-discretisation of the heat equation problem with $N$ internal grid points and  the control on the right boundary:
  \begin{equation}
    \label{heat}
      \left\{         \begin{split}
     \partial_{t} v -  \nu \partial_{xx} v&=0, \qquad (t, x) \in \oi 0T \times \oi 01,\\
     v(t, 0)&=0, \quad v(t, 1)=u_\nu(t),\\
      v(0, x)&=v_0.\\
     \end{split}\right.
     \end{equation}
     The parameter $\nu$ represents the diffusion coefficient and is supposed to range within the set $\lN=\zi{1}2$.  The system satisfies the Kalman's rank condition for any $\nu\in \lN$ and any target time $T$.
     
We aim to control the system  from the initial state $v_0(x)=\sin(\pi x)$ to zero in time $T=0.1$.

      The greedy control algorithm has been  applied for the system of dimension $N=50$ with  $\eps=0.0001$, and the uniform discretisation of $\lN$ in $k=100$ values.
  
  The algorithm stops after only three iterations, choosing  parameter values (out of 100 eligibles) in the following order:
      $$
   1.00 ,\,  1.18,\,   1.45\,.
   $$
        The elapsed time is 213 seconds (obtained on the same machine used for example 7.1).

   The corresponding three minimisers $\ph_i^0$ are used for  constructing approximate controls  for all  parameter values   by means of formulas \eqref{approx_control} and \eqref{approx_min}. Here we present the results obtained for $\nu=\sqrt 2$.
 \begin{figure}
      \vskip -4cm  
      \centerline{  \includegraphics[scale=0.5]{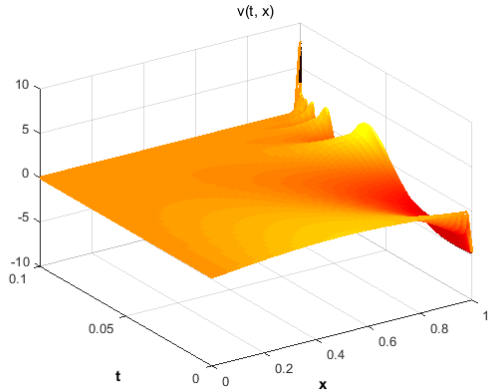}}
       \vskip -4cm
   \caption{\label{heat-3D}Evolution of the solution to the semi-discretised problem  \eqref{heat} governed by the approximate control $u_\nu^\star$ for $\nu=\sqrt 2$.}
  \end{figure}

\begin{figure}
      \vskip -5cm  
      \centerline{ \includegraphics[scale=0.5]{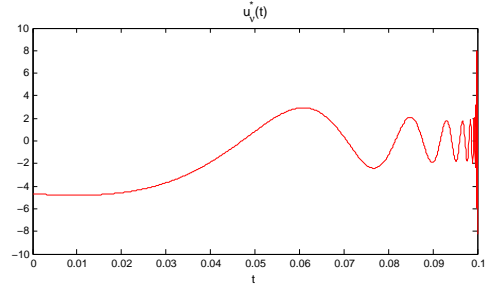}}
      \vskip -5cm
       \caption{\label{heat-control} Evolution of
     the approximate control for 
     $\nu=\sqrt 2$.}
  \end{figure}

 Evolution of the solution, presented by a 3-D plot is given by   Figure \ref{heat-3D}. The system is driven to the zero state within the error $|x(T)|=1\cdot 10^{-5}$. The corresponding approximate control profile is depicted on Figure \ref{heat-control}, exhibiting rather strong oscillations when  approaching the target time -- an intrinsic feature for the heat equation. 

    Let us mention that the elapsed time for constructing the approximate control $u_\nu^\star$ was 1.5 seconds only, compared to 37 seconds required to construct the optimal exact control. This huge difference is due to a rather small number of snapshots required to obtain approximation performance within the given accuracy,  and it  reconfirms computational efficiency of the greedy approach when  compared to the standard one.

   The choice of a rather tiny value of the target time $T=0.1$ is due to the strong dissipation effect of the heat equation, providing an exponential   solution decay even in the absence of any control. For the same reason, the  algorithm  stopped after only three iterations although the precision was set rather high with  $\eps=0.0001$ (for the wave problem 24 iterations were required in order to produce uniform approximation within the error 0.5).

\section{Conclusion and perspectives}

\begin{itemize}

 \item  {\it Exact controllability.} In our analysis the exact control problem is relaxed to the approximate one by letting a tolerance $\eps$ on the target at time $T$, what is realistic in applications. Note however that this approximation of the final controllability condition is achieved by ensuring a sufficiently small error on the approximation of exact controls. Thus, in practice, the methods we have built are of direct application to the problem of exact controllability as well as the tolerance $\eps$ tends to zero.
   
 \item  {\it Comparison of greedy and naive approach.} As we have seen the weak greedy algorithm we have developed has various advantages with respect to the naive one, the latter  consisting of taking a fine enough uniform mesh on  the set of eligible parameters $\lN$  and calculating all corresponding controls. 
   Although the greedy  approach requires an (expensive) search for distinguished parameter values, in the finite-dimensional case its   computational cost is smaller than for the naive  one.  Meanwhile,  in the infinite-dimensional context it leads to algorithms with an optimal convergence rate.

 \item   {\it Model reduction:} The computational cost of the greedy methods, as developed here, remains high. This is so because on the search of the new distinguished parameter values, to compute the surrogate, we are obliged to solve fully the adjoint system first and then the controlled one. It would be interesting to build cheaper (from a computational viewpoint)  surrogates. This would require implementing model reduction techniques. But this, as far as we know, has non been successfully done in the context of controllability. Model reduction methods have been implemented for elliptic and parabolic optimal control problems (see \cite{GNVPL}, \cite{GK}, and \cite{KG}) but, as far as we know, they have not been developed in the context of controllability.
   
 \item {\it Parameter dependence.}
Most  existing practical realisations of  (weak) greedy algorithms assume an  affine parameter dependence  presented at the end of Section \ref{infinite} (\cite{CD15,D15,KG}). The assumption provides a cheaper computation of a surrogate and reduces the cost of the search for a next snapshot. 

By the contrary,   the greedy control algorithm presented in this paper allows for very general parameter  dependence, still providing an efficient algorithm beating the naive approach. Smooth dependence is needed in order to obtain sharp convergence rates. Still the algorithms can be applied for models depending on the parameters in a rough manner.

   The greedy control algorithm and corresponding approximation results of Section 4, as well as solution to our Problem is derived under assumption of Lipschitz continuity of the control system entries with respect to the parameter.
   However, as explained in Remark \ref{cont_dep}, all these results remain
   valid if  we merely assume continuous dependence with respect to the parameter. The Lipschitz regularity is chosen just  to  simplify the presentation
   and to obtain explicit expressions  by means  of the given approximation error $\eps$ and
   the Lipschitz constants.
   
   A generalisation of the approach to  problems  involving discontinuous parameter dependence would require an approximation by piecewise smooth functions and  separate  applications of the algorithm on each subinterval fulfilling the  continuity assumption.

   The additional analyticity requirement is imposed in Section \ref{infinite} in order to deduce convergence analysis of the algorithm. 
   Namely,  transfer of Kolmogorov width from the set of parameters into the set of solutions (or controls) is provided so far only by Theorem \ref{greedy_rates2}, which requires analyticity assumption, as well as polynomial decay of Kolmogorov widths.
    Once this result is expanded to a more general setting, it can be applied directly to the presented greedy control algorithm, and provide convergence estimates under weaker    regularity assumptions, without modifying the very algorithm.

   However, the present lack of an analogous result under more general  assumptions does not prevent possible applications of the greedy control algorithm to non-analytic data.
   It can still be implemented, practically at no-risk. Namely, in the worst case when dimension of approximation space reaches number of eligible parameters,  one will calculate controls for all these parameter values, which corresponds to the  naive approach. And this will cost the same as applying naive method directly from the beginning.
But there are many reasons to believe the algorithm  will stop before this ultimate point.

\item {\it Waves versus heat:} In our numerical experiments we have observed that the greedy algorithm is more efficient for the semi-discrete heat equation than  for the wave one, in the sense that for the first one it stops after $n=3$ iterations while for the second one it goes up to $n=24$. This is a natural result to be expected. Indeed, for the heat equation, because of the intrinsic dissipativity of the system, the high frequency components play a minor role when, as in our experiments, dealing with approximate controllability. Thus, even if the algebraic dimension of both systems, waves and heat, are the same, in practice the relevant dimensions for the heat equation is much smaller, thus explaining the faster convergence of the greedy algorithm for heat like equations.

\item {\it Extension of greedy control to PDE systems.} As we have seen the methods and ideas developed in this paper can be easily adapted to infinite-dimensional ODEs. The same ideas and methods can be implemented without major changes  on the  PDE setting. However, as we have seen above, in order to quantify the convergence rate of greedy algorithms, the analytic dependence of the controls with respect to parameters plays a key role. As far as we know, this issue has not been thoroughly addressed in the literature, i.e. that of whether or not controls for a given PDE controllability problem depend analytically on the parameters entering in the system. And at this level, the type of PDE under consideration can play a major role. Indeed,  for heat equations, even when the unknown parameters enter in the principal part of the diffusion operator, the analytic dependence of the controls can be expected. This is not the case for wave equations (see \cite{AK}). Indeed, note that even for the simplest first order transport equation solutions fail to depend analytically on the velocity of propagation in a natural $L^2$ or Sobolev energy setting. This issue requires significant further investigation.

 \item  {\it Robust control.} The analysis and methods developed in this paper apply to control problems where the data to be controlled are prescribed a priori, either depending on the parameters or not.
 In that sense, our results are the analogue to what has been done for the greedy approximation of PDE solutions with given data (right hand side terms and boundary values) as in  \cite{BCDDPW},  \cite{CDV} and \cite{CD15}, among others. 
 
  From the point of view of applications it would be interesting to develop greedy methods for control 
of potential application for all possible data to be controlled. This would require to establish a strategy for identifying the most relevant snapshots of the parameters $\nu$ for the approximation of the Gramians $\Lambda(\nu)$, within the space of bounded linear operators. And this, in turn, requires identifying efficient surrogates. This is an interesting problem to be addressed.

Actually, as far as we know, this has not been done even in the context of the solvability of elliptic problems. The work done so far, as mentioned above, addresses approximation issues for given specific data. But the problem of applying greedy algorithms to approximate the resolvent operators, in a robust manner and independently of   the data entering in the PDE, has not been addressed.

\end{itemize}

\noindent {\bf Acknowledgements.} The authors acknowledges A. Cohen and R. DeVore  for fruitful discussions, and J. Loh\' eac for valuable support with the numerical simulations.

\end{document}